\documentclass[12pt]{amsart}
\usepackage[T1]{fontenc}
\usepackage[utf8]{inputenc}
\usepackage[english]{babel}
\usepackage{amssymb,amsmath}
\usepackage{suffix}
\usepackage{graphicx}
\usepackage{multirow}
\usepackage{geometry}
\usepackage{natbib}
\usepackage{hyperref}
\usepackage{pdfcomment}
\newcommand{\bbR}{\mathbb{R}}
\newcommand{\calC}{\mathcal{C}}
\newcommand{\calP}{\mathcal{P}}
\newcommand{\ext}{\mathop{\mathrm{ext}}}

\newcommand{\diff}{\mathrm{d}}
\newcommand{\dt}{\diff t}
\newcommand{\dd}[2][{}]{\frac{\diff#1}{\diff#2}}
\newcommand{\ddt}{\dd[]t}
\newcommand{\pp}[2][{}]{\frac{\partial#1}{\partial#2}}
\WithSuffix\newcommand\pp*[2][{}]{\partial#1/\partial#2}
\newcommand{\dq}{\dot{q}}
\newcommand{\leg}{\mathop{\mathrm{leg}}\nolimits}
\newcommand{\cf}{\emph{cf.}}
\newcommand{\eg}{\emph{e.g.}}

\newcommand{\quand}{\quad\text{and}\quad}
\newcommand{\qquand}{\qquad\text{and}\qquad}

\title[HOVI, a polynomial approach]{Higher Order Variational Integrators:\\a polynomial approach}
\author[C. M. Campos]{Cédric M. Campos}
\address{Instituto de Ciencias Matemáticas\\
         Campus de Cantoblanco, UAM\\
         Calle Nicolás Cabrera, 15\\
         28049 Madrid (Spain)}
\thanks{\emph{Subject.} On the occasion of the 23rd Congress on Differential Equations and Applications, CEDYA 2013, held in Castellón, Spain.\\
        \indent \emph{Acknowledgments.} This work has been partially supported by MEC (Spain) Grants MTM2010-21186-C02-01, MTM2011-15725E, the ICMAT Severo Ochoa project SEV-2011-0087 and the European project IRSES-project ``Geomech-246981''. Besides, the author wants to specially thank professors Sina Ober-Blöbaum, from the University of Paderborn, and David Martín de Diego, from the Instituto de Ciencias Matemáticas, for fruitful conversations on the subject and their constant support.}

\hypersetup{
  linkcolor=blue,
  citecolor=blue,
  filecolor=blue,
  urlcolor=MidnightBlue,
  pdftitle={Higher Order Variational Integrators: a polynomial approach},
  pdfauthor={C. M. Campos},
  pdfsubject={XXIII Congreso de Ecuaciones Diferenciales y Aplicaciones, CEDYA 2013},
  pdfkeywords={geometric mechanics, Euler-Lagrange equation, discrete mechanics, geometric integration, variational integrator, higher-order, Runge-Kutta methods},
  pdfcreator={GNU Emacs 24.3.1},
  pagebackref=true,
}

\begin{document}

\maketitle

\begin{abstract}
We reconsider the variational derivation of symplectic partitioned Runge-Kutta schemes. Such type of variational integrators are of great importance since they integrate mechanical systems with high order accuracy while preserving the structural properties of these systems, like the symplectic form, the evolution of the momentum maps or the energy behaviour. Also they are easily applicable to optimal control problems based on mechanical systems as proposed in \cite{ObJuMa10}.

Following the same approach, we develop a family of variational integrators to which we refer as symplectic Galerkin schemes in contrast to symplectic partitioned Runge-Kutta. These two families of integrators are, in principle and by construction, different one from the other. Furthermore, the symplectic Galerkin family can as easily be applied in optimal control problems, for which \cite{CaJuOb12} is a particular case.
\end{abstract}

\section{Introduction}\label{sec:intro}
In recent years, much effort in designing numerical methods for the time integration of (ordinary) differential equations has been put into schemes which are \emph{structure preserving} in the sense that important \emph{qualitative} features of the original dynamics are preserved in its time discretization, \cf\ the recent monograph \cite{HaLuWa02}. A particularly elegant way to, \eg\ derive symplectic integrators, is by discretizing Hamilton's principle as suggested by \cite{Su90,Ve88}, see also \cite{MaWe01,SaCa94}.

However most part of the theory and examples rely on second order schemes, hence some effort must still be put into the development of accurate higher order schemes that, in long term simulations, can drastically reduce the overall computational cost. A clear example are the so called symplectic partitioned Runge-Kutta methods that integrate mechanical systems driven by a Lagrangian $L\colon(q,\dq)\in\bbR^n\times\bbR^n\mapsto L(q,\dq)\in\bbR$ and, possibly, by a force $f\colon(q,\dq)\in\bbR^n\times\bbR^n\mapsto p=f(q,\dq)\in\bbR^n$. A detailed study of such methods can be found in \cite{HaLuWa02}.

A partitioned Runge-Kutta method is an $s$-stage integrator $(q_0,p_0)\mapsto(q_1,p_1)$ given by the equations
\begin{subequations} \label{eq:pRK}
\begin{align}
q_1 =& q_0 + h\sum_{j=1}^sb_j\dot Q_j\,, & p_1 =& p_0 + h\sum_{j=1}^s\bar b_j\dot P_j\,,\\
Q_i =& q_0 + h\sum_{j=1}^sa_{ij}\dot Q_j\,, & P_i =& p_0 + h\sum_{j=1}^s\bar a_{ij}\dot P_j\,,\\
P_i =& \pp[L]{\dot q}(Q_i,\dot Q_i)\,, & \dot P_i =& \pp[L]{\dot q}(Q_i,\dot Q_i) + f(Q_i,\dot Q_i)\,,
\end{align}
\end{subequations}
where $(b_j,a_{ij})$ and $(\bar b_j,\bar a_{ij})$ are two different Runge-Kutta methods associated to collocation points $0\leq c_1<\ldots<c_s\leq1$ and time step $h>0$.

It is shown that the previous integrator is \emph{symplectic} whenever the two sets of coefficients satisfy the relations
\begin{subequations} \label{eq:spRK.coeff}
\begin{eqnarray}
b_i\bar a_{ij}+\bar b_ja_{ji} &=& b_i\bar b_j \,,\\
b_i = \bar b_i \,.
\end{eqnarray}
\end{subequations}
In fact, it can be derived as a geometric variational integrator, see \cite{MaWe01,Su90}.

The paper is structured as follows: Section~\ref{sec:dm} is a short introduction to Discrete Mechanics and Section~\ref{sec:hovi} describes the variational derivation of higher order schemes using polynomial collocation. Finally we briefly enumerate the relations and differences between symplectic partitioned Runge-Kutta schemes and symplectic Galerkin ones in Section~\ref{sec:rel} and conclude by outlining future research directions in Section~\ref{sec:conclusion}.

\section{Discrete Mechanics and Variational Integrators}\label{sec:dm}

One of the main subjects of Geometric Mechanics is the study of dynamical systems governed by a Lagrangian. Typically they consider mechanical systems with \emph{configuration manifold} $Q$ together with a \emph{Lagrangian function} $L\colon TQ\to\bbR$, where the associated \emph{state space} $TQ$ describes the position and velocity of a particle moving in the system. Usually, the Lagrangian takes the form of kinetic minus potential energy, $ L(q,\dq) = K(q,\dq) - V(q) = \frac12\,\dq^T\cdot M(q)\cdot\dq - V(q)$, for some (positive definite) \emph{mass matrix} $M(q)$.

A consequence of the \emph{principle of least action}, also known as \emph{Hamilton's principle}, establishes that the natural motions $q\colon[0,T]\to Q$ of the system are characterized by the celebrated \emph{Euler-Lagrange equation} (refer to \cite{AbMa78}),
\begin{equation} \label{eq:EL-equation}
\pp[L]{q} - \ddt\pp[L]{\dq} = 0 \,.
\end{equation}
When the Lagrangian is \emph{regular}, that is when the velocity Hessian matrix $\partial^2 L/\partial \dq^2$ is non-degenerate, the Lagrangian induces a well defined map, the \emph{Lagrangian flow}, $F_L^t\colon TQ\to TQ$ by $F_L(q_0,\dq_0) := (q(t),\dq(t))$, where $q\in\calC^2([0,T],Q)$ is the unique solution of the Euler-Lagrange equation \eqref{eq:EL-equation} with initial condition $(q_0,\dq_0)\in TQ$. By means of the \emph{Legendre transform} $\leg_L:(q,\dq)\in TQ\mapsto (q,p=\pp[L]{\dq}|_{(q,\dq)})\in T^*Q$, where $T^*Q$ is the \emph{phase space} of positions plus momenta, one may transform the Lagrangian flow into the \emph{Hamiltonian flow} $F_H^t\colon T^*Q\to T^*Q$ defined by $F_H^t(q_0,p_0) := \leg_L(q(t),\dq(t))$.

Moreover, different preservation laws are present in these systems. For instance the Hamiltonian flow preserves the natural symplectic structure of $T^*Q$ and the total energy of the system. Also, if the Lagrangian possess Lie group symmetries, then \emph{Noether's theorem} asserts that some quantities are conserved, like for instance the linear momentum and/or the angular momentum.

Discrete Mechanics is, roughly speaking, a discretization of Geometric Mechanics theory. As a result, one obtains a set of discrete equations equivalent to the Euler-Lagrange equation \eqref{eq:EL-equation} above but, instead of a direct discretization of the ODE, the latter are derived from a discretization of the base objects of the theory, the state space $TQ$, the Lagrangian $L$, etc. In fact, one seeks for a sequence $\{(t_0,q_0),(t_1,q_1),\dots,\allowbreak(t_n,q_n)\}$ that approximates the actual trajectory $q(t)$ of the system ($q_k\approx q(t_k)$), for a constant time-step $h=t_{k+1}-t_k>0$.

A \emph{variational integrator} is an iterative rule that outputs this sequence and it is derived in an analogous manner to the continuous framework. Given a discrete Lagrangian $L_d\colon Q\times Q\to\bbR$, which is in principle thought to approximate the continuous Lagrangian action over a short time
\[ L_d(q_k,q_{k+1}) \approx \int_{t_k}^{t_{k+1}}L(q(t),\dq(t))\dt \,, \]
one applies a variational principle to derive the well-known discrete Euler-Lagrange (DEL) equation,
\begin{equation} \label{eq:DEL}
D_1L_d(q_k,q_{k+1}) + D_2L_d(q_{k-1},q_k) = 0 \,,
\end{equation}
where $D_i$ stands for the partial derivative with respect to the $i$-th component. The equation defines an integration rule of the type $(q_{k-1},q_k) \mapsto (q_k,q_{k+1})$, however if we define the pre- and post-momenta
\begin{equation} \label{eq:discrete.legendre.transform}
p_k^- := -D_1L_d(q_k,q_{k+1})  \qquand  p_k^+ := D_2L_d(q_{k-1},q_k) \,,
\end{equation}
the Euler-Lagrange equation \eqref{eq:EL-equation} is read as the momentum matching $p_k^-=p_k^+=:p_k$ and defines an integration rule of the type $(q_k,p_k)\mapsto(q_{k+1},p_{k+1})$.

The nice part of the story is that the integrators derived in this way naturally preserve (or nearly preserve) the quantities that are preserved in the continuous framework, the symplectic form, the total energy and, in presence of symmetries, the linear and/or angular momentum (for more details, see \cite{MaWe01}). Furthermore, other aspects of the continuous theory can be ``easily'' adapted, symmetry reduction \cite{CaCeDi12,CoJiMa12,IgMaMa08}, constraints \cite{JoMu09,KoMaSu10}, control forces \cite{CaJuOb12,ObJuMa10}, etc.

\section{Higher Order Variational Integrators}\label{sec:hovi}

Higher order variational integrators for time dependent or independent systems (HOVIt) are a class of integrators that, by using a multi-stage approach, aim at a higher order accuracy on the computation of the natural trajectories of a mechanical system while preserving some intrinsic properties of such systems. In particular, symplectic-partitioned Runge-Kutta methods (spRK) and, what we call here, symplectic Galerkin methods (sG) are $s$-stage variational integrators of order up to $2s$.

The derivation of these methods follows a general scheme. For a fixed time step $h$, one considers a series of points $q_k$, refereed as macro-nodes. Between each couple of macro-nodes $(q_k,q_{k+1})$, one also considers a set of micro-data, the $s$ stages: For the particular cases of sG and spRK methods, micro-nodes $Q_1,\dots,Q_s$ and micro-velocities $\dot Q_1,\dots,\dot Q_s$, respectively. Both macro-nodes and micro-data (micro-nodes or micro-velocities) are required to satisfy a variational principle, giving rise to a set of equations, which properly combined, define the final integrator.

In what follows, we will use the following notation: Let $0\leq c_1<\ldots<c_s\leq1$ denote a set of collocation times and consider the associated Lagrange polynomials and nodal weights, that is,
\[ l^j(t) := \prod_{i\neq j}\frac{t-c_i}{c_j-c_i} \qquand b_j := \int_0^1l^j(t)\dt \,, \]
respectively. Note that the pair of $(c_i,b_i)$'s define a quadrature rule and that, for appropriate $c_i$'s, this rule may be a Gaussian-like quadrature, for instance, Gauss-Legendre, Gauss-Lobatto, Radau or Chebyschev.

Now, for the sake of simplicity and independently on the method, we will use the same notation for the nodal coefficients. We define for spRK and sG, respectively,
\[ a_{ij} := \int_0^{c_i}l^j(t)\dt  \qquand  a_{ij} := \dd[l^j]t\Big|_{c_i} \,. \]
Moreover, for spRK, we will also use the nodal weights and coefficients $(\bar b_j,\bar a_{ij})$ given by Equation \eqref{eq:spRK.coeff} and, for sG, the source and target coefficients
\[ \alpha^j := l^j(0)  \qquand  \beta^j:=l^j(1) \,. \]

Finally, if $L$ denotes a Lagrangian from $\bbR^n\times\bbR^n$ to $\bbR$, then we define
\[ P_i := \pp[L]{\dq}\Big|_i = \pp[L]{\dq}\Big|_{(Q_i,\dot Q_i)}  \qquand  \dot P_i := \pp[L]{q}\Big|_i = \pp[L]{q}\Big|_{(Q_i,\dot Q_i)} \,, \]
where $(Q_i,\dot Q_i)$ are couples of micro-nodes and micro-velocities given by each method. Besides, $D_i$ will stand for the partial derivative with respect to the $i$-th component.

\subsection{Symplectic-Partitioned Runge-Kutta Methods}\label{sec:sprk}
Although the variational derivation of spRK methods in the framework of Geometric Mechanics is an already known fact (see \cite{MaWe01} for an ``intrinsic'' derivation, as the current, or \cite{HaLuWa02} for a ``constrained'' one), we present it here again in order to ease the comprehension of and the comparison with sG methods below.

Given a point $q_0\in\bbR^n$ and vectors $\{\dot Q_i\}_{i=1,\dots,s}\subset\bbR^n$, we define the polynomial curves
\[ \dot Q(t) := \sum_{j=1}^sl^j(t/h)\dot Q_j \qquand Q(t) := q_0+h\sum_{j=1}^s\int_0^{t/h}l^j(\tau)\diff\tau\dot Q_j \,. \]
We have
\[ \dot Q_i = \dot Q(h\cdot c_i) \quand Q_i := Q(h\cdot c_i) = q_0 + h\sum_{j=1}^sa_{ij}\dot Q_j \,. \]
Note that the polynomial curve $Q$ is uniquely determined by $q_0$ and $\{\dot Q_i\}_{i=1,\dots,s}$. In fact, it is the unique polynomial curve $Q$ of degree $s$ such that $Q(0)=q_0$ and $\dot Q(h\cdot c_i)=\dot Q_i$. However, if we define the configuration point
\begin{equation} \label{eq:spRK.target.node}
q_1 := Q(h\cdot1) = q_0 + h\sum_{j=1}^sb_j\dot Q_j
\end{equation}
and consider it fixed, then $Q$ is uniquely determined by $q_0$, $q_1$ and the $\dot Q_i$'s but one. Namely, take any $1\leq i_0 \leq s$ such that $b_{i_0}\neq0$ and fix it, then
\[ \dot Q_{i_0} = \left(\frac{q_1-q_0}h-\sum_{j\neq i_0}b_j\dot Q_j\right)/b_{i_0} \]
and we obtain the following relations that will be useful in what comes next
\[ \pp[(Q_i,Q_{i_0})]{(q_0,Q_j,q_1)} =
   \left( \begin{array}{ccc}
       0 & \delta_i^j & 0\\
     -\frac1{hb_{i_0}} & -\frac{b_j}{b_{i_0}} & \frac1{hb_{i_0}}
   \end{array} \right) \,,\ i,j\neq i_0. \]


We now define the multi-vector discrete Lagrangian
\[ L_d(\dot Q_1,\dots,\dot Q_s) := h\sum_{i=1}^sb_iL(Q_i,\dot Q_i) \,. \]
Although not explicitly stated, it also depends on $q_0$. The two-point discrete Lagrangian is then
\[ L_d(q_0,q_1) := \ext_{\calP^s([0,h],\bbR^n,q_0,q_1)}L_d(\dot Q_1,\dots,\dot Q_s) \]
where $\calP^s([0,h],\bbR^n,q_0,q_1)$ is the space of polynomials $Q$ of order $s$ from $[0,1]$ to $\bbR^n$ such that $Q(0)=q_0$ and $Q(h)=q_1$ and the vectors $\dot Q_i$'s determine such polynomials as discussed above. The extremal is realized by a polynomial $Q\in\calP^s([0,h],\bbR^n,q_0,q_1)$ such that
\begin{equation} \label{eq:DEL.micro.spRK}
\delta L_d(\dot Q_1,\dots,\dot Q_s)\cdot(\delta \dot Q_1,\dots,\delta \dot Q_s) = 0
\end{equation}
for any variations $(\delta \dot Q_1,\dots,\delta \dot Q_s)$, taking into account that
\[ \delta q_0=\delta q_1=0 \qquand \delta \dot Q_{i_0} = \sum_{j\neq i_0}\pp[\dot Q_{i_0}]{\dot Q_j}\delta \dot Q_j \,. \]
For convenience, the previous equation is developed afterwards.

By the momenta-matching rule \eqref{eq:discrete.legendre.transform}, we have that
\begin{eqnarray*}
-p_0 &=& D_1L_d(q_0,q_1) \\
     &=& D_{i_0}L_d(\dot Q_1,\dots,\dot Q_s)\pp*[\dot Q_{i_0}]{q_0} + D_{q_0}L_d(\dot Q_1,\dots,\dot Q_s)\\
     &=& -D_{i_0}L_d(\dot Q_1,\dots,\dot Q_s)/(hb_{i_0}) + D_{q_0}L_d(\dot Q_1,\dots,\dot Q_s) \,,\\
     & &\\
 p_1 &=& D_2L_d(q_0,q_1) \\
     &=& D_{i_0}L_d(\dot Q_1,\dots,\dot Q_s)\pp*[\dot Q_{i_0}]{q_1}\\
     &=& D_{i_0}L_d(\dot Q_1,\dots,\dot Q_s)/(hb_{i_0}) \,.
\end{eqnarray*}
where $D_{q_0}$ stands for the partial derivative with respect to $q_0$. Combining both equations, we obtain that
\[ D_{i_0}L_d(\dot Q_1,\dots,\dot Q_s) = hb_{i_0}p_1 \qquand  p_1 = p_0 + D_{q_0}L_d(\dot Q_1,\dots,\dot Q_s) \,. \]
Coming back to Equation \eqref{eq:DEL.micro.spRK}, we have that
\begin{multline*}
\delta L_d(\dot Q_1,\dots,\dot Q_s)\cdot(\delta \dot Q_1,\dots,\delta \dot Q_s) = \\
\begin{split}
&= \sum_{j\neq i_0}D_jL_d(\dot Q_1,\dots,\dot Q_s)\delta \dot Q_j + D_{i_0}L_d(\dot Q_1,\dots,\dot Q_s)\delta \dot Q_{i_0}\\
&= \sum_{j\neq i_0}\left(D_jL_d(\dot Q_1,\dots,\dot Q_s) + \pp[\dot Q_{i_0}]{\dot Q_j}D_{i_0}L_d(\dot Q_1,\dots,\dot Q_s)\right)\delta \dot Q_j \,.
\end{split}
\end{multline*}
Therefore, for $j\neq i_0$, we have that
\[ D_jL_d(\dot Q_1,\dots,\dot Q_s) = b_j/b_{i_0}\cdot D_{i_0}L_d(\dot Q_1,\dots,\dot Q_s)  \,. \]
Thus, for any $j=1,\dots,s$, we have that
\begin{equation} \label{eq:spRK.DEL}
D_jL_d(\dot Q_1,\dots,\dot Q_s) = hb_jp_1 \,.
\end{equation}

The integrator is defined by
\begin{subequations}
\begin{gather}
D_jL_d(\dot Q_1,\dots,\dot Q_s) = hb_jp_1 \,, \\
q_1 = q_0 + h\sum_{j=1}^sb_j\dot Q_j \,, \\
p_1 = p_0 + D_{q_0}L_d(\dot Q_1,\dots,\dot Q_s) \,.
\end{gather}
\end{subequations}
Besides, using the definition of the discrete Lagrangian, we have
\begin{eqnarray*}
D_jL_d(\dot Q_1,\dots,\dot Q_s)
&=& h\sum_{i=1}^sb_i\left( \pp[L]{q}\Big|_i\pp[Q_i]{\dot Q_j} + \pp[L]{\dq}\Big|_i\pp[\dot Q_i]{\dot Q_j} \right) \\
&=& h^2\sum_{i=1}^sb_ia_{ij}\dot P_i + hb_jP_j \,, \\
D_{q_0}L_d(\dot Q_1,\dots,\dot Q_s)
&=& h\sum_{i=1}^sb_i\left( \pp[L]{q}\Big|_i\pp[Q_i]{q_0} + \pp[L]{\dq}\Big|_i\pp[\dot Q_i]{q_0} \right) \\
&=& h\sum_{i=1}^sb_i\dot P_i \,.
\end{eqnarray*}
Therefore, we may write
\begin{gather*}
P_j = p_0 + h\sum_{i=1}^sb_i(1-a_{ij}/b_j)\dot P_i = p_0 + h\sum_{i=1}^s\bar a_{ji}\dot P_i \,, \\
p_1 = p_0 + h\sum_{i=1}^sb_i\dot P_i = p_0 + h\sum_{i=1}^s\bar b_i\dot P_i \,,
\end{gather*}
were $\bar a_{ij}$ and $\bar b_i$ are given by Equation \eqref{eq:spRK.coeff}.

In summary, we have recovered all the equations that define the spRK integrator without forces, Equation \eqref{eq:pRK} and \eqref{eq:spRK.coeff}, that is
\begin{subequations} \label{eq:spRK}
\begin{align}
q_1 =& q_0 + h\sum_{j=1}^sb_j\dot Q_j\,, & p_1 =& p_0 + h\sum_{j=1}^s\bar b_j\dot P_j\,,\\
Q_i =& q_0 + h\sum_{j=1}^sa_{ij}\dot Q_j\,, & P_i =& p_0 + h\sum_{j=1}^s\bar a_{ij}\dot P_j\,,\\
P_i =& \pp[L]{\dot q}(Q_i,\dot Q_i)\,, & \dot P_i =& \pp[L]{\dot q}(Q_i,\dot Q_i) \,.
\end{align}
\end{subequations}

\subsection{Symplectic Galerkin Methods (of 0th kind)}\label{sec:sg}
Galerkin methods are a class of methods to transform a problem given by a continuous operator (such as a differential operator) to a discrete problem. As such, spRK methods falls into the scope of this technique and could be also classified as ``symplectic Galerkin'' methods (of 1st kind). However, we want to stress here the difference between what is called spRK and what we here refer as sG. The wording should not be confused by the one used in \cite{MaWe01}.

Given points $\{Q_i\}_{i=1,\dots,s}\subset\bbR^n$, we define the polynomial curves
\[ Q(t) := \sum_{j=1}^sl^j(t/h)Q_j \qquand \dot Q(t) := \frac1h\sum_{j=1}^s\dot l^j(t/h)Q_j \,. \]
We have
\[ Q_i = Q(h\cdot c_i) \quand \dot Q_i := \dot Q(h\cdot c_i) = \frac1h\sum_{j=1}^sa_{ij}Q_j \,. \]
Note that the polynomial curve $Q$ is uniquely determined by the points $\{Q_i\}_{i=1,\dots,s}$. In fact, it is the unique polynomial curve $Q$ of degree $s$ such that $Q(h\cdot c_i)=Q_i$. However, if we define the configuration points
\begin{equation} \label{sG.source-target.nodes}
q_0 := Q(h\cdot0) = \sum_{j=1}^s\alpha^jQ_j  \qquand  q_1 := Q(h\cdot1) = \sum_{j=1}^s\beta^jQ_j
\end{equation}
and consider them fixed, then $Q$ is uniquely determined by $q_0$, $q_1$ and the $Q_i$'s but a couple. For instance, we may consider $Q_1$ and $Q_s$ as functions of the others, since the relations \eqref{sG.source-target.nodes} define a system of linear equations where the coefficient matrix has determinant $\gamma:=\alpha^1\beta^s-\alpha^s\beta^1\neq0$ (if and only if $c_1\neq c_s$). More precisely,
\[
   \left(
   \begin{array}{c}
     Q_1\\
     Q_s
   \end{array}
   \right)
   =
   \frac1\gamma
   \left(
   \begin{array}{cc}
      \beta^s & -\alpha^s\\
     -\beta^1 &  \alpha^1
   \end{array}
   \right)
   \left(
   \begin{array}{c}
     q_0-\sum_{j=2}^{s-1}\alpha^jQ_j\\
     q_1-\sum_{j=2}^{s-1}\beta^jQ_j
   \end{array}
   \right) \,,
\]
from where we obtain the following relations that will be useful in what comes next
\[ \pp[(Q_1,Q_i,Q_s)]{(q_0,Q_j,q_1)} =
   \frac1\gamma
   \left( \begin{array}{ccc}
      \beta^s & \alpha^s\beta^j-\alpha^j\beta^s & -\alpha^s\\
       0 & \gamma\delta_i^j & 0\\
     -\beta^1 & \alpha^j\beta^1-\alpha^1\beta^j &  \alpha^1
   \end{array} \right) \,,\ i,j=2,\dots,s-1. \]
We will also take into account that
\[ \pp[\dot Q_i]{(q_0,Q_j,q_1)} = \frac1h\sum_{k=1}^sa_{ik}\pp[Q_k]{(q_0,Q_j,q_1)} \,,\ i=1,\dots,s,\ j=2,\dots,s-1.\]

We now define the multi-point discrete Lagrangian
\[ L_d(Q_1,\dots,Q_s) := h\sum_{i=1}^sb_iL(Q_i,\dot Q_i) \,. \]
The two-point discrete Lagrangian is then
\[ L_d(q_0,q_1) := \ext_{\calP^s([0,h],\bbR^n,q_0,q_1)}L_d(Q_1,\dots,Q_s) \]
where $\calP^s([0,h],\bbR^n,q_0,q_1)$ is the space of polynomials $Q$ of order $s$ from $[0,1]$ to $\bbR^n$ such that the points $Q_i$'s determine such polynomials as discussed above. The extremal is realized by a polynomial $Q\in\calP^s([0,h],\bbR^n,q_0,q_1)$ such that
\begin{equation} \label{eq:DEL.micro.sG}
\delta L_d(Q_1,\dots,Q_s)\cdot(\delta Q_1,\dots,\delta Q_s) = 0
\end{equation}
for any variations $(\delta Q_1,\dots,\delta Q_s)$, taking into account that
\[ \delta q_0=\delta q_1=0  \qquand \delta Q_i = \sum_{j=2}^{s-1}\pp[Q_i]{Q_j}\delta Q_j \,,\ i=1,s. \]
For convenience, the previous equation is developed afterwards.

By the momenta-matching rule \eqref{eq:discrete.legendre.transform}, we have that
\begin{eqnarray*}
-p_0 &=& D_1L_d(q_0,q_1) \\
     &=& \sum_{j=1}^sD_jL_d(Q_1,\dots,Q_s)\pp[Q_j]{q_0}\\
     &=& \beta^s/\gamma \cdot D_1L_d(Q_1,\dots,Q_s) - \beta^1/\gamma \cdot D_sL_d(Q_1,\dots,Q_s)\\
     & &\\
 p_1 &=& D_2L_d(q_0,q_1) \\
     &=& \sum_{j=1}^sD_jL_d(Q_1,\dots,Q_s)\pp[Q_j]{q_1}\\
     &=& -\alpha^s/\gamma \cdot D_1L_d(Q_1,\dots,Q_s) + \alpha^1/\gamma \cdot D_sL_d(Q_1,\dots,Q_s)
\end{eqnarray*}
By a linear transformation of both equations, we obtain
\begin{align*}
  D_1L_d(Q_1,\dots,Q_s) &= -\alpha^1p_0 + \beta^1p_1 \ \ \text{and}\\
  D_sL_d(Q_1,\dots,Q_s) &= -\alpha^sp_0 + \beta^sp_1 \,.
\end{align*}
Coming back to Equation \eqref{eq:DEL.micro.sG}, we have that
\begin{multline*}
\delta L_d(Q_1,\dots,Q_s)\cdot(\delta Q_1,\dots,\delta Q_s) = \\
\begin{split}
&= D_1L_d(Q_1,\dots,Q_s)\delta Q_1 + \sum_{j=2}^{s-1}D_jL_d(Q_1,\dots,Q_s)\delta Q_j + D_sL_d(Q_1,\dots,Q_s)\delta Q_s\\
&= \sum_{j=2}^{s-1}\left[D_1L_d(Q_1,\dots,Q_s)\pp[Q_1]{Q_j} + D_jL_d(Q_1,\dots,Q_s) + D_sL_d(Q_1,\dots,Q_s)\pp[Q_s]{Q_j}\right]\delta Q_j
\end{split}
\end{multline*}
Therefore, for $j=2,\dots,s-1$, we obtain
\begin{eqnarray*}
\gamma D_jL_d
&=& (\alpha^j\beta^s-\alpha^s\beta^j)D_1L_d + (\alpha^1\beta^j-\alpha^j\beta^1)D_sL_d\\
&=& (\alpha^j\beta^s-\alpha^s\beta^j)(-\alpha^1p_0+\beta^1p_1) + (\alpha^1\beta^j-\alpha^j\beta^1)(-\alpha^sp_0+\beta^sp_1)\\
&=& (\alpha^1\beta^s-\alpha^s\beta^1)(-\alpha^jp_0+\beta^jp_1) \,.
\end{eqnarray*}
Thus, for any $j=1,\dots,s$, we have that
\begin{equation} \label{eq:sG.DEL}
D_jL_d(Q_1,\dots,Q_s) = -\alpha^jp_0 + \beta^jp_1 \,.
\end{equation}

The integrator is defined by
\begin{subequations}
\begin{gather}
D_jL_d(Q_1,\dots,Q_s) = -\alpha^jp_0 + \beta^jp_1 \,,\ j=1,\dots,s;\\
q_0 = \sum_{j=1}^s\alpha^jQ_j  \qquand q_1 = \sum_{j=1}^s\beta^jQ_j
\end{gather}
\end{subequations}
Besides, using the definition of the discrete Lagrangian, we have
\begin{eqnarray*}
D_jL_d(Q_1,\dots,Q_s)
&=& h\sum_{i=1}^sb_i\left( \pp[L]{q}\Big|_i\pp[Q_i]{\dot Q_j} + \pp[L]{\dq}\Big|_i\pp[\dot Q_i]{\dot Q_j} \right) \\
&=& hb_j\dot P_j + \sum_{j=1}^sb_ia_{ij}P_i \,.
\end{eqnarray*}
Therefore, we may simply write
\[ hb_j\dot P_j + \sum_{j=1}^sb_ia_{ij}P_i = -\alpha^jp_0 + \beta^jp_1 \,. \]

In summary and for a proper comparison, we write the equations that define the sG integrator (without forces) in a pRK fashion, that is
\begin{subequations} \label{eq:sG}
\begin{align}
q_0 =& \sum_{j=1}^s\alpha^jQ_j\,, & q_1 =& \sum_{j=1}^s\beta^jQ_j\,,\\
Q_i =& \frac1h\sum_{j=1}^sa_{ij}Q_j\,, & \dot P_i =& \frac{\beta^ip_1-\alpha^ip_0}{h\bar b_i} + \frac1h\sum_{j=1}^s\bar a_{ij}P_j\,,\\
P_i =& \pp[L]{\dq}(Q_i,\dot Q_i)\,, & \dot P_i =& \pp[L]{\dq}(Q_i,\dot Q_i)\,,
\end{align}
\end{subequations}
where $b_ia_{ij}+\bar b_j\bar a_{ji}=0$ and $b_i=\bar b_j$.

We remark that Equation \eqref{eq:sG.DEL} generalizes the ones obtained in \cite{CaJuOb12,Le04}, where the collocation times are chosen such that $c_1=0$ and $c_s=1$, which is a rather particular case.

\section{Relations between spRK and sG}\label{sec:rel}
First of all, it is worth to say that, with a little bit of extra technicalities, one can easily include forces into both schemes. As a result, one would only need to redefine in \eqref{eq:spRK} and \eqref{eq:sG}
\[ \dot P_i = \pp[L]{\dq}(Q_i,\dot Q_i) + f(Q_i,\dot Q_i) \,, \]
where $f\colon(q,\dq)\in\bbR^n\times\bbR^n\mapsto p=f(q,\dq)\in\bbR^n$ is the external force.

Both methods can be considered of Galerkin type. In this sense, spRK and sG could be refereed as a symplectic Galerkin integrators of 1st and 0th kind, respectively, since spRK is derived from the 1st derivative of an extremal polynomial and sG from the polynomial itself. At this point, a very natural question could arise: Actually are spRK and sG two different integrator schemes? Even though the derivations of both methods are quite similar, they are in general different (although they could coincide for particular choices of the Lagrangian, the collocation times and the integral quadrature). A weak but still fair argument for this is that, at each step, spRK relies on the determination of the micro-velocities $\dot Q_i$, while sG does so on the micro-nodes $Q_i$. All the other ``variables'' are then computed form the determined micro-data.

With respect to the accuracy of the schemes, for any Gaussian quadrature (Gauss-Legendre, Gauss-Lobatto, Radau and Chebyschev) and any method (spRK and sG), the schemes have convergence order $2s-2$, expect for the combination of Gauss-Lobatto together with spRK which is $2s$, being $s$ the number of internal stages.

\begin{table}[h]
  \centering
  \begin{tabular}{cc|c|c|}
    \cline{3-4}
    & & spRK & sG \\
    \cline{1-4}
    \multicolumn{2}{|c|}{micro-data} & $\dot Q_i$ & $Q_i$ \\
    \cline{1-4}
    \multicolumn{2}{|c|}{polynomial degree} & $s$ & $s-1$ \\
    \cline{1-4}
    \multicolumn{2}{|c|}{variational eq.'s} & $s+1$ & $s$ \\
    \cline{1-4}
    \multicolumn{2}{|c|}{extra equations} & $1$ & $2$ \\
    \cline{1-4}
    \multicolumn{1}{|c|}{\multirow{4}{*}{quadrature}} & Gauss-Legendre & $2s$ & $2s-2$ \\
    \cline{2-4}
    \multicolumn{1}{|c|}{\multirow{4}{*}{}} & Gauss-Lobatto & $2s-2$ & $2s-2$ \\
    \cline{2-4}
    \multicolumn{1}{|c|}{\multirow{4}{*}{}} & Radau & $2s-2$ & $2s-2$ \\
    \cline{2-4}
    \multicolumn{1}{|c|}{\multirow{4}{*}{}} & Chebyschev & $2s-2$ & $2s-2$ \\
    \cline{1-4}
    & & \multicolumn{2}{|c|}{order method} \\
    \cline{3-4}
  \end{tabular}

  \caption{Comparison of $s$-stage variational integrators.}
  \label{tab:comparison}
\end{table}

Let's finish underlying that sG, as spRK, is inherently symplectic.

\section{Conclusions}\label{sec:conclusion}
In this work, by revisiting the variational derivation of spRK methods \cite{HaLuWa02,MaWe01}, we have presented a new class of higher order variational integrators within the family of Galerkin schemes. These integrators are symplectic per construction and, therefore, well suited for long term simulations, where the higher order accuracy of the schemes can be exploited to reduce the overall computational cost. Also, they can be easily adapted to implement constraints or symmetry reduction \cite{CaCeDi12,IgMaMa08,MaWe01} or, together with an NLP solver, to integrate optimal control problems \cite{CaJuOb12,ObJuMa10}.

For the future, the sG schemes deserves a proper analysis to establish the actual differences with respect to spRK schemes and results the convergence rates. And to further take advantage of the high accuracy of the methods, we envisage to design time adaptive algorithms. Joint work with O. Junge (TUM), S. Ober-Blöbaum (UPB) and E. Trélat (CNRS-UPMC) on the control direction has already started in order to generalize \cite{CaJuOb12} and clarify some aspects of \cite{Ha00}. As well, we have begun digging along the lines of constrained systems and higher order Lagrangians.

\bibliographystyle{plainnat}
\bibliography{hovi}

\end{document}